\documentclass[11pt]{amsart}

\usepackage{mathptmx} 
\usepackage[scaled=0.90]{helvet} 
\usepackage{courier} 
\normalfont
\usepackage[T1]{fontenc}
\usepackage{amsmath,amsthm,amssymb}

\usepackage[curve,tips]{xypic}
\SelectTips{eu}{11}
\UseTips

\newtheorem{theorem}{Theorem}[section]
\newtheorem{proposition}[theorem]{Proposition}

\newtheorem{lemma}[theorem]{Lemma}

\newtheorem*{theoremA}{Theorem A}
\newtheorem*{theoremB}{Theorem B}
\newtheorem*{theoremC}{Theorem C}
\newtheorem*{corollaryA}{Corollary A}

\theoremstyle{definition}
\newtheorem{definition}[theorem]{Definition}
\newtheorem{remark}[theorem]{Remark}

\renewcommand{\setminus}{\smallsetminus}
\renewcommand{\emptyset}{\varnothing}

\newcommand{\normal}{\lhd}
\newcommand{\mono}{\rightarrowtail}
\newcommand{\epi}{\twoheadrightarrow}
\newcommand{\iso}{\cong}
\newcommand{\Q}{{\mathbb Q}}
\newcommand{\Z}{{\mathbb Z}}
\newcommand{\F}{{\mathbb F}}
\newcommand{\fpinfty}{\operatorname{FP}_\infty}
\newcommand{\fp}{\operatorname{FP}}
\newcommand{\cd}{\operatorname{cd}}

\newcommand{\PD}{\operatorname{PD}}
\newcommand{\Tor}{\operatorname{Tor}}
\newcommand{\Gal}{\operatorname{Gal}}
\newcommand{\colim}{{\displaystyle\lim_{\buildrel\longrightarrow\over\lambda}\ }}
\newcommand{\colimf}%
{{\displaystyle\lim_{\buildrel\longrightarrow\over{H\in\mathcal S}}\ }}
\newcommand{\limf}%
{{\displaystyle\lim_{\buildrel\longleftarrow\over{H\in\mathcal S}}\ }}
\newcommand{\invlim}{{\displaystyle\lim_{\longleftarrow}}}

\newcommand{\Ider}{\operatorname{Ider}}
\newcommand{\Der}{\operatorname{Der}}
\newcommand{\Comm}{\operatorname{Comm}}
\newcommand{\Mo}{\operatorname{Mod}}
\newcommand{\HYPHEN}{\operatorname{-}}
\newcommand{\Mod}{\operatorname{\Mo\HYPHEN}}

\title[A spectral sequence and applications]%
{A generalization of the\\ Lyndon--Hochschild--Serre spectral sequence\\ with applications to group cohomology\\ and decompositions of groups}
\author{P. H. Kropholler}
\address{University of Glasgow, Glasgow G12 8QQ}
\email{p.h.kropholler@maths.gla.ac.uk}

\begin{document}

\begin{abstract}
We set up a Grothendieck spectral sequence which generalizes the Lyndon--Hochschild--Serre spectral sequence for a group extension $K\mono G\epi Q$ by allowing the normal subgroup $K$ to be replaced by a subgroup, or  family of subgroups which satisfy a weaker condition than normality. This
is applied to establish a decomposition theorem for certain groups as fundamental groups of graphs of Poincar\'e duality groups. We further illustrate the method by proving a cohomological vanishing theorem which applies for example to Thompson's group $F$.
\end{abstract}

\maketitle


For motivation consider a finitely generated group $G$ of cohomological dimension $n+1$ over $\Z$, where $n$ is a non-negative integer. Suppose that $H$ is a polycyclic subgroup of $G$ which has Hirsch length $n$. If $H$ is normal in $G$ then a Lyndon--Hochschild--Serre spectral sequence corner argument shows that the quotient group $G/H$ has more than one end. By a classical result of Stallings it follows that $G/H$ splits as a non-trivial amalgamated free product or HNN-extension over a finite subgroup and so $G$ itself splits over a subgroup in which $H$ has finite index. If $H$ is not normal then there may be no splitting over any subgroup related to $H$. In \cite{ds} Dunwoody and Swenson consider a situation like this, but rather than imposing a cohomological dimension condition they consider the notion of {\em codimension one} by which they mean that the number of ends of the pair $(G,H)$ is $\ge2$. Here, we can think of the spectral sequence corner argument as a device which proves that when $H$ is normal then it has codimension one in the sense of Dunwoody and Swenson. 
One could define a new notion of {\em cohomological codimension one} to mean that $\cd H=\cd G-1$ so that the LHS corner argument may be summarized by saying that if $H$ is normal
\begin{center}
cohomological codimension one $\implies$ codimension one.
\end{center}
It is at once clear that the notion of cohomological codimension one is altogether much weaker than the notion defined in terms of ends. Dunwoody and Swenson obtain substantial and very complete splitting results under a codimension one hypothesis with a polycyclic-by-finite subgroup. On the other hand, cohomological codimension one is a commonplace. For example if $G$ has cohomological dimension $2$ and $H$ is any infinite cyclic subgroup then we have cohomological codimension one but associated splittings would be regarded as exceptional rather than typical.

One of the main results in this paper is a splitting theorem which is proved by showing that codimension one can be deduced from cohomological codimension one in a little more generality than the case of a normal subgroup. What we do is to construct a Grothendieck spectral sequence generalizing the LHS spectral sequence but where normality of the subgroup is replaced by a weaker condition. The result, Theorem B, is a new splitting theorem which should prove useful in the general study of JSJ decompositions for groups as discussed in \cite{ss} by Scott and Swarup. Our basic idea is to design the new spectral sequence so that the corner argument will work  in the same way as it did in the classical case.

We set the scene with some definitions before describing the spectral sequence.
Let $G$ be a group. Two subgroups $H$ and $K$ of $G$ are said to be {\em commensurable} if and only if $H\cap K$ has finite index in both $H$ and $K$. The {\em commensurator} $$\Comm_G(H)$$ of a subgroup $H$ is defined to be the set of $g\in G$ such that $H$ and $H^g$ are commensurable. This is a subgroup of $G$ which contains the normalizer $N_G(H)$ of $H$. 
\begin{definition}\label{def:near normal}
We say that $H$ is {\em near-normal} in $G$ if and only if $\Comm_G(H)=G$. 
\end{definition}
Near-normal subgroups arise naturally: for example, $SL_n(\Z)$ is a near-normal subgroup of $SL_n(\Q)$.
We shall also use the following terminology.
\begin{definition}\label{def:admissible}
Let $G$ be a group. By an {\em admissible family} of subgroups of $G$ we shall mean a family $\mathcal S$ which is closed under conjugation and such that the intersection of any finite number of members of $\mathcal S$ contains a member of $\mathcal S$.
\end{definition}
The simplest instance of an admissible family is the family $\{K\}$ consisting of a single normal subgroup $K$. For a more interesting example, suppose that $H$ is a near-normal subgroup of $G$ and take $\mathcal S$ to be the family of all subgroups commensurable with $H$. 

Our first step is a Grothendieck spectral sequence which generalizes the Lyndon--Hochschild--Serre spectral sequence for a group extension $$K\mono G\epi Q$$ by allowing the normal subgroup $K$ to be replaced by a near-normal subgroup, or more generally, any admissible family of subgroups. 
Here are two key examples which the reader may wish to keep in mind while reading the paper.
\begin{itemize}
\item
The Baumslag--Solitar group
$G=\langle x,y:\ y^{-1}x^2y=x^3\rangle$ and $K=\langle x\rangle$. 
\newline
We replace $K$ by the family of all subgroups which contain a positive power of $x$.
\item
Thompson's group
$F=\langle x_0^{\phantom{1}},x_1^{\phantom{1}},x_2^{\phantom{1}},\dots:\ x_i^{-1}x_jx_i^{\phantom{1}}=x_{j+1}^{\phantom{1}}\ (i<j)\rangle$, and

$K=\langle x_{1}^{\phantom{1}}x_{0}^{-1},x_{3}^{\phantom{1}}x_{2}^{-1},x_{5}^{\phantom{1}}x_{4}^{-1},x_{7}^{\phantom{1}}x_{6}^{-1},\dots\rangle$.
\newline 
We replace $K$ by the family of all finite intersections of conjugates of $K$ by words of even exponent sum in the generators $x_i$. 
\end{itemize}

We introduce the basic abelian categories and establish our new spectral sequence in \S1 and further study cohomological issues in \S2. Our major application to the theory of duality groups and decompositions of groups is presented in \S3. We discuss some of the implications of this for Poincar\'e duality in \S4. Another very simple application of the spectral sequence method is included in the very short section \S5, where we prove a cohomological vanishing theorem: this result provides a simple and elegant way of seeing why $H^*(F,\Z F)=0$ for Thompson's group $F$. We conclude the paper with a discussion of the connection between our methods and those used in Galois cohomology.

I would like to acknowledge with thanks several useful conversations with Martin Dunwoody and Gadde Swarup regarding the application in \S3 and with John Wilson and Pavel Zalesskii regarding the closing remarks in \S6.
Both these and the referee's advice have led to improvements and to the correction of some flaws in the original draft.

\section{The spectral sequence}

By a Grothendieck spectral sequence we mean a composite functor spectral sequence as described by Cartan and Eilenberg (\cite{cartaneilenberg}, Chapter XVII \S7). The Lyndon--Hochschild--Serre spectral sequence (\cite{cartaneilenberg}, Chapter XVI \S6 (6)) can be viewed this way: see for example Rotman's account, Theorem 11.45 of \cite{rotman}.
Consider a group extension $K\mono G\epi Q$. Write $\Mod\Z G$ for the category of right $\Z G$-modules. 
The zeroth cohomology functor $H^0(G,{\phantom M})$ is the fixed point functor $({\phantom M})^G$ and the LHS spectral sequence arises from the factorization of this functor through the category of $\Z Q$-modules as illustrated below:
\[
\xymatrix{
\Mod\Z G\ar[rr]^{({\phantom M})^G}\ar[dr]_{H^0(K,{\phantom M})=({\phantom M})^K}&&\Mod\Z.\\
&\Mod\Z Q\ar[ur]_{H^0(Q,{\phantom M})=({\phantom M})^Q}\\
}
\]
One needs to know that the module categories are abelian categories with enough injectives and that the $K$-fixed point functor carries injective $G$-modules to injective $Q$-modules.

Now suppose that $\mathcal S$ is an admissible family of subgroups of $G$. 
\begin{definition}\label{new cat}
We write $$\Mod\Z G/\mathcal S$$ for the full subcategory of $\Mod\Z G$ whose objects are those $\Z G$-modules $M$ for which
$M=\bigcup_{H\in\mathcal S}M^H.$ 
\end{definition}

\begin{remark}\label{rem} Although this definition makes perfect sense for any family $\mathcal{S}$ of subgroups, it behaves particularly well when $\mathcal{S}$ is admissible and we shall always assume that this is so. Given this assumption and an arbitrary $\Z G$-module $M$ then $\bigcup_{H\in\mathcal S}M^H$ is a $\Z G$-submodule of $M$ which belongs to the new category. We write $H^0(\mathcal{S},M)$ for this submodule. Notice that the definition of admissible family is designed exactly so that this works: $H^0(\mathcal{S},M)$ inherits an action of $G$ because $\mathcal{S}$ is closed under conjugation and $H^0(\mathcal{S},M)$ is an additive subgroup of $M$ because $\mathcal{S}$ is downwardly directed. Thus $H^0(\mathcal{S},M)$ is an object of the new category $\Mod\Z G/\mathcal S$ and {\em by definition} all objects of $\Mod\Z G/\mathcal S$ arise this way.
\end{remark}
It is easy to see that $\Mod\Z G/\mathcal S$ is an abelian category:
it will be the intermediary for our spectral sequence, replacing $\Mod\Z Q$.
Note that in case $\mathcal S$ consists of a single (necessarily normal) subgroup $K$ then $\Mod\Z G/\mathcal S$ is naturally equivalent to the category of right modules for the quotient group $Q=G/K$. 

We now have two functors. The first is mentioned already in Remark \ref{rem}:
the assignment
$$H^0(\mathcal S,{\phantom M}):\Mod\Z G\to\Mod\Z G/\mathcal S$$ 
defined by
$$H^0(\mathcal S,M)=\bigcup_{H\in\mathcal S}M^H$$
is functorial in $M$. Secondly, we can restrict the $G$-fixed point functor to the new category so we have a functor
$$H^0(G/\mathcal S,{\phantom M}):\Mod\Z G/\mathcal S\to\Mod\Z$$
defined by
$$H^0(G/\mathcal S,M)=M^G.$$
\begin{remark}\label{rem2}
The $G$-fixed point functor $({\phantom M})^G:\Mod\Z G\to\Mod\Z$ now factors through 
$\Mod\Z G/\mathcal S$ as the composite: $H^0(\mathcal S,{\phantom M})$ followed by
$H^0(G/\mathcal S,{\phantom M})$
as illustrated in the diagram below.
\[
\xymatrix{
\Mod\Z G\ar[rr]^{({\phantom M})^G}\ar[dr]_{H^0(\mathcal S,{\phantom M})}&&\Mod\Z.\\
&\Mod\Z G/\mathcal S\ar[ur]_{H^0(G/\mathcal S,{\phantom M})}\\
}
\]
This will provide the basis for a Grothendieck spectral sequence.
\end{remark}
Beware that the notation $H^0(G/\mathcal S,{\phantom M})$ is not intended to imply the construction of any kind of object $G/\mathcal S$, but is simply notation for the functor. This functor necessarily has to be distinguished from $H^0(G,{\phantom M})$ which has a different domain. The notation {\em is} intended to suggest an analogy with the classical situation when $\mathcal S$ consists of a single normal subgroup $K$.
The analogy works well, and raises the interesting question whether there is any kind of natural object which deserves to be named $G/\mathcal S$. In \S6 we show that a certain completion of $G$ appears to be the object one should expect.
\begin{lemma}\label{injectives}\ 
\begin{enumerate}
\item
If $I$ is an injective $\Z G$-module then $H^0(\mathcal S,I)$ is injective in $\Mod\Z G/\mathcal S$.
\item
$\Mod\Z G/\mathcal S$ has enough injectives.
\item
The functor $H^0(\mathcal S,{\phantom M})$ has right derived functors 
$H^n(\mathcal S,{\phantom M})$ and there are natural isomorphisms
$$H^n(\mathcal S,{\phantom M})\iso\colimf H^n(H,M).$$
\end{enumerate}
\end{lemma}
\begin{proof} (i) Let $I$ be an injective $\Z G$-module. To show that $H^0(\mathcal S,I)$ is injective we need to address the extension problem as illustrated below:
\[
\xymatrix{
0\ar[r]&M\ar[d]_\phi\ar[r]&N\ar@{-->}[dl]^{\widehat\phi{\text?}}\\
&H^0(\mathcal S,I)&\\
}
\]
Using the injectivity of $I$ we can find $\widehat\phi$ to make a commutative diagram
\[
\xymatrix{
0\ar[r]&M\ar[d]_\phi\ar[r]&N\ar[ddl]^{\widehat\phi}\\
&H^0(\mathcal S,I)\ar[d]&\\
&I&
}
\]
and since $N$ belongs to the subcategory, $\widehat\phi$ has image in $H^0(\mathcal S,I)$. This proves (i). For part (ii), observe that we can embed any $\Z G/\mathcal S$-module $M$ into an injective $\Z G$-module and then we have $M=H^0(\mathcal S,M)\hookrightarrow H^0(\mathcal S,I)$: by part (i), we have now embedded $M$ into an injective object of $\Mod\Z G/\mathcal S$.

(iii) Derived functors are defined in the standard way using injective resolutions 
over $\Z G$, applying the functor $H^0(\mathcal S,{\phantom M})$ and passing to the cohomology of the resulting cochain complex.
The isomorphism is easily established.
\end{proof}
\begin{lemma}\label{lem:der quot}
The functor $H^0(G/\mathcal S,{\phantom M})$ has right derived functors 
$H^n(G/\mathcal S,{\phantom M})$.
\end{lemma}
\begin{proof}
This time we work with injective resolutions in the category $\Mod\Z G/\mathcal S$,
apply the $G$-fixed point functor and pass to cohomology. In general, there is no simple interpretation of the derived functors.
\end{proof}
Remark \ref{rem2}, Lemma \ref{injectives} and Lemma \ref{lem:der quot} together provide the ingredients necessary for a Grothendieck spectral sequence and our main tool is established:
\begin{theoremA}\label{A}
Let $G$ be a group and let $\mathcal S$ be an admissible family of subgroups. 
There is a Grothendieck spectral sequence 
$$H^p(G/\mathcal S,H^q(\mathcal S,M))\implies H^{p+q}(G,M)$$
which is natural in the $G$-module $M$.
\end{theoremA}
As a routine feature of any first/third quadrant spectral sequence we have the following: (see for example \cite{rotman}, Theorem 11.43).
\begin{corollaryA}\label{corollaryA}
With $G$ and $\mathcal S$ as above, 
\begin{enumerate}
\item there is a five term exact sequence analogous to the standard inflation-restriction sequence:
{\small{$$0\to H^1(G/\mathcal S,H^0(\mathcal S,M))\to H^1(G,M)\to
H^1(\mathcal S,M)^G\to H^2(G/\mathcal S,H^0(\mathcal S,M))\to H^2(G,M);$$}}
\item the inflation map $$H^n(G/\mathcal S,M)\to H^n(G,M)$$ (which is defined for $M$ in the subcategory $\Mod\Z G/\mathcal S$) is an isomorphism when $n=0$ and is injective when $n=1$.
\end{enumerate}
\end{corollaryA}

\section{Continuity of functors}

We need to consider continuity issues for the new functors. 
\begin{definition}\label{def:cont}
Let $\mathcal C$ and $\mathcal D$ be abelian categories with all small filtered colimits.
Let $F:\mathcal C\to\mathcal D$ be a functor. We say that $F$ is {\em continuous at zero} if and only if
$$\colim F(M_\lambda)=0$$ whenever $(M_\lambda)$ is a small filtered colimit system in $\mathcal C$ which is {\em vanishing}: i.e.
$$\colim M_\lambda=0.$$
\end{definition}
First we recall a basic result from Bieri's notes, essentially the content of  (\cite{bieri-qmw}, Theorem 1.3 (i)$\iff$(iiib)):
\begin{lemma}\label{bieri1}
A group $G$ is of type $\fp_n$ if and only if the cohomology functors $H^i(G,{\phantom M})$ are continuous at zero for all $i\le n$.
\end{lemma}
\begin{lemma}\label{lem:fpinfty}
Let $\mathcal S$ be an admissible family in $G$. If all members of $\mathcal S$ have type $\fpinfty$ then the functors $H^n(\mathcal S,{\phantom M})$ are continuous at zero for all $n$.
\end{lemma}
\begin{proof} By Lemma \ref{bieri1}, the $\fpinfty$ condition guarantees that the functors $H^{m}(H,{\phantom M})$ are continuous at zero for all $H\in\mathcal S$ and all $m$.
The result now follows from the natural isomorphism of Lemma \ref{injectives}(iii).
\end{proof}

The next lemma is a version of Strebel's criterion \cite{strebel}. 
\begin{lemma}\label{strebel}
Let $G$ be a group of finite cohomological dimension. Suppose that for all vanishing filtered colimit systems $(P_\lambda)$ of projective modules $P_\lambda$, and all $m\in\Z$,
$$\colim H^m(G,P_\lambda)=0.$$ Then $G$ is of type $\fp$.
\end{lemma}
\begin{proof}
We shall use the following notation: write $FM$ for the free module on the underlying set of non-zero elements of a module $M$. Then $F$ is functorial in $M$ and the inclusion $M\setminus\{0\}\hookrightarrow M$ induces a natural surjection $FM\epi M$ whose kernel, $\Omega M$, is also functorial. Moreover both $F$ and $\Omega$ take vanishing filtered colimit systems to vanishing filtered colimit systems.

By Lemma \ref{bieri1} it suffices to prove that 
$H^m(G,{\phantom M})$ is continuous at zero for all $m$. This is proved by downward induction on $m$: it is trivial for all $m$ greater than the cohomological dimension of $G$, so we fix $m$ and assume as inductive hypothesis that $H^{m+1}(G,{\phantom M})$ is continuous at zero. Let $(M_\lambda)$ be a vanishing filtered colimit system of modules.
Then we have short exact sequences
$$\Omega M_\lambda\mono FM_\lambda\epi M_\lambda.$$ On passing to cohomology and taking colimits we have the exact sequence
$$\colim H^m(G,FM_\lambda)\to \colim H^m(G,M_\lambda)\to \colim H^{m+1}(G,\Omega M_\lambda).$$
We need to prove that the central group here is zero. The right hand group vanishes by induction and the left hand group vanishes by hypothesis. The result follows from exactness.
\end{proof}

\begin{lemma}\label{lem:van1}
Let $L$ be a near-normal subgroup of type $\fpinfty$ and infinite index in a group $G$. 
Let $M$ be a $\Z L$-module. Let $\mathcal S$ be the set of subgroups commensurable with $L$. Then 
$$\left(H^m(\mathcal S,M\otimes_{\Z L}\Z G)\right)^G=0$$ for all integers $m$.
\end{lemma}
\begin{proof}
We proceed in steps.
\begin{enumerate}
\item[Step 1.] The case $m=0$.
\end{enumerate}
We use only on the fact that all members of $\mathcal S$ have infinite index in $G$. 
Note that the calculation simplifies because
$$\left(H^0(\mathcal S,M\otimes_{\Z L}\Z G)\right)^G=(M\otimes_{\Z L}\Z G)^G,$$
the set of $G$-fixed points.
Let $T$ be a right transversal to $L$ in $G$, i.e. $G$ is the disjoint union of the cosets $Lt$, for $t\in T$. Any non-zero element of $M\otimes_{\Z L}\Z G$ has a unique expression as a finite sum 
$$m_1\otimes t_1+\dots+m_s\otimes t_s$$
where the $m_i\in M$ are non-zero and the $t_i$ are distinct elements of $T$.
Let $X=Lt_1\cup\dots\cup Lt_s$. We can choose $g\in G$ such that $Xg\cap X=\emptyset$. To see this, suppose for a contradiction that $Xg\cap X=\emptyset$ for all $g$. Then $G$ is the union of the sets $t_i^{-1}Lt_j$ over all $i,j$: this expresses $G$ as a finite union of cosets of subgroups and implies that at least one of the subgroups has finite index in $G$ which is contrary to our assumption. For a $g$ such that $Xg\cap X=\emptyset$ we have
$$(m_1\otimes t_1+\dots+m_s\otimes t_s)g\ne m_1\otimes t_1+\dots+m_s\otimes t_s$$
so there are no non-zero fixed points.
\begin{enumerate}
\item[Step 2.] In case $M$ is injective as a $\Z L$-module and $m\ge1$.
\end{enumerate}
We prove the stronger statement that for any $H$ in $\mathcal S$,
$$H^m(H,M\otimes_{\Z L}\Z G)=0.$$
This argument uses both the commensurability and the type $\fpinfty$.
Fix any $H$. Mackey decomposition yields
$$M\otimes_{\Z L}\Z G\iso\bigoplus_tMt\otimes_{\Z[L^t\cap H]}\Z H$$
as $\Z H$-modules, where $t$ runs over a set of $(L,H)$ double coset representatives in $G$. Each summand $Mt\otimes_{\Z[L^t\cap H]}\Z H$ is injective as an $H$-module because $M$ is injective over $L$ and, using commensurability, 
$L^t\cap H$ has finite index in $H$. In positive dimensions, the cohomology of any group vanishes on injective modules. Here we can take advantage of the fact that $H$ has type $\fpinfty$ to see that the cohomology also vanishes on arbitrary direct sums of injective modules.
\begin{enumerate}
\item[Step 3.] The general case.
\end{enumerate}
The general case can be deduced by dimension shifting. We use induction on $m$, the case $m=0$ being covered by step (i). Choose any short exact sequence
$$M\mono I \epi M'$$ with $I$ injective over $L$. Then we have a short exact sequence of induced modules:
$$M\otimes_{\Z L}\Z G\mono I\otimes_{\Z L}\Z G\epi M'\otimes_{\Z L}\Z G.$$ Using the long exact sequence of cohomology together with step (ii) reduces the $m$-dimensional matter for $M$ to the $(m-1)$-dimensional matter for $M'$ and the result follows by induction. Of course, in case $m=1$, Step 2 does not cover everything but then Step 1 can be used as well.
\end{proof}

\begin{proposition}\label{bieri2}
Let $n$ be natural number.
Let $G$ be a finitely generated group of cohomological dimension $\le n+1$. Let $K$ be near-normal subgroup of infinite index in $G$ such that:
\begin{enumerate}
\item
$K$ is of type $\fp$;
\item
$H^m(K,P)=0$ whenever $P$ is a projective module and $m\ne n$.
\end{enumerate}
Then $G$ is of type $\fp$.
\end{proposition}
\begin{proof} Let $\mathcal S$ be the family of all subgroups commensurable with $K$. Let $P$ be a projective $\Z G$-module. All the members of $\mathcal S$ inherit properties (i) and (ii) and it follows from Lemma \ref{injectives}(iii) and Lemma \ref{bieri1} that $\mathcal S$ itself inherits the properties as well:
\begin{enumerate}
\item
$H^m(\mathcal S,{\phantom M})$ is continuous at zero for all $m$;
\item
$H^m(\mathcal S,P)=0$ whenever $P$ is a projective $\Z G$-module and $m\ne n$.
\end{enumerate}
The spectral sequence of Theorem A therefore collapses to a single column and we find that 
$$H^{n+1}(G,P)\iso H^1(G/\mathcal S,H^n(\mathcal S,P)),$$
$$H^{n}(G,P)\iso H^0(G/\mathcal S,H^n(\mathcal S,P)),$$
and 
$$H^{m}(G,P)=0$$ when $m\notin\{n,n+1\}$ because of (ii) for $m<n$ and the constraint on dimension of $G$ for $m>n+1$.
By definition, $H^0(G/\mathcal S,H^n(\mathcal S,{\phantom M}))=(H^n(\mathcal S,{\phantom M}))^G$ and this functor vanishes on all induced modules (i.e. modules of the form $A\otimes\Z G$ where $A$ is an abelian group) by Lemma \ref{lem:van1}. Since projective modules are direct summands of free modules which are in turn examples of induced modules, we have
$$H^{n}(G,P)=0.$$
Now we can apply Strebel's criterion Lemma \ref{strebel}. If $(P_\lambda)$ is a vanishing filtered colimit system of projective modules then we only have to check that
$$\colim H^{n+1}(G,P_\lambda)=0$$
and this will follow from the isomorphism 
$$\colim H^{n+1}(G,P_\lambda)=\colim H^1(G/\mathcal S,H^n(\mathcal S,P_\lambda))$$
together with the observation that both the functors $H^1(G/\mathcal S,{\phantom M})$ and $H^n(\mathcal S,{\phantom M})$ are continuous at zero. The second observation is part of (i) above. The first observation follows from the injectivity of the inflation maps
$$H^1(G/\mathcal S,H^n(\mathcal S,P_\lambda))\to 
H^1(G,H^n(\mathcal S,P_\lambda))$$
at the start of the five term exact sequence, Corollary A: note that $H^1(G,{\phantom M})$ is continuous at zero because $G$ is finitely generated.
\end{proof}

\section{A Decomposition Theorem}

The goal of this section is to prove a decomposition theorem for certain groups. We take this in two stages. The first stage, Theorem B, is easy to state. The second stage Theorem C requires a little more introduction although it is easy to prove by combining Theorem B with results \cite{bf} of Bestvina and Feighn.

\begin{theoremB}\label{B} Let $n$ be a fixed natural number. Suppose that $G$ is a group with the following properties:
\begin{itemize}
\item[$(\alpha)$]
$G$ is finitely generated;
\item[$(\beta)$]
$G$ has cohomological dimension $\le n+1$;
\item[$(\gamma)$]
$G$ has a  near-normal $\PD^n$-subgroup $K$.
\end{itemize}
Then either $K$ has finite index in $G$ or $G$ splits over a subgroup commensurable with $K$. 
\end{theoremB}

We recall some definitions. An {\em $n$-dimensional duality group} (over $\Z$) is a group $G$ which affords a dualizing module $D$ so that there are natural isomorphisms
$$H^i(G,M)\iso \Tor^{\Z G}_{n-i}(M,D)$$ for all $i$.
{\em Poincar\'e duality groups} ($\PD^n$-groups) are duality groups for which $D$ has underlying additive group $\Z$. We say that a group $G$ {\em splits} over a subgroup $H$ if and only if $G$ is isomorphic to a free product with amalgamation $K*_HL$ with $K\ne H\ne L$ or and HNN-extension $K*_H$. According to the standard theory of group actions on trees which is described by Serre \cite{serre} and Dicks--Dunwoody \cite{dicks}, $G$ splits over $H$ if and only if there is a $G$-tree with no fixed vertex, one orbit of edges and in which $H$ is the stabilizer of one of the edges.

The proof of Theorem B uses Dunwoody's graph cutting methods, \cite{dunwoody}.

\begin{proof}[Proof of Theorem B]\ 
We may as well assume that $K$ has infinite index in $G$.
Let $\mathcal S$ be the family of all subgroups commensurable with $K$. Then $\mathcal S$ is an admissible family of $\PD^n$-subgroups of $G$ all of which have infinite index in $G$.
The hypotheses of Proposition \ref{bieri2} are satisfied and therefore $G$ has type $\fp$.
The next step is to prove 
\begin{itemize}
\item[{\bf Claim 1.}] $G$ is an $(n+1)$-dimensional duality group over any field $k$.
\end{itemize}
To be an $(n+1)$-dimensional duality group over a (non-zero) commutative ring $k$ it is necessary and sufficient that the following three conditions hold.
\begin{itemize}
\item
$G$ is of type $\fp$ over $k$ and of cohomological dimension $\le n+1$.
\item
The cohomology groups $H^i(G,kG)$ are zero for $i\le n$.
\item $H^{n+1}(G,kG)$ is flat as a $k$-module.
\end{itemize}
When these conditions hold, $G$ has dualizing module $H^{n+1}(G,kG)$: this is a left $kG$-module via the left action of $G$ on $kG$.
Since $kG$ is an instance of an induced $\Z G$-module, the first two conditions are already established. The third condition is automatically satisfied if $k$ is a field. We now work over the field $\F$ of two elements and show how to identify $H^n(\mathcal S,\F G)$ with a certain set of subsets of $G$.

Let $\mathcal P$ denote the powerset of $G$. This is viewed as an $(\F G,\F G)$-bimodule with symmetric difference of subsets providing the additive structure and left/right multiplication by elements of $G$ for the action.
Let $\mathcal B$ be the set of subsets $B$ of $G$ which satisfy the following condition
\begin{itemize}
\item
There is a subgroup $H\in\mathcal S$ and a finite subset $F$ of $G$ such that $B=HFH$. I.e. $B$ is a finite union of double cosets of some $\mathcal S$-subgroup.
\end{itemize}
The set $\mathcal B$ is an $(\F G,\F G)$-sub-bimodule.
\begin{itemize}
\item[{\bf Claim 2.}]
$H^n(\mathcal S,\F G)$ is isomorphic to $\mathcal B$ as a $(\F G,\F G)$-bimodule.
\end{itemize}
We need to understand the connecting maps which are involved in the colimit formulation for $H^*(\mathcal S,\F G)$. When $H\subseteq L$ are $\PD^n$-groups then  $H$ has finite index in $L$ and there are commutative diagrams
\[
\xymatrix{
H^i(L,M)\ar[r]\ar[d]^{\text{Res}}&\Tor_{n-i}^{\F L}(M,\F)\ar[d]^{\text{Tr}}\\
H^i(H,M)\ar[r]&\Tor_{n-i}^{\F H}(M,\F)\\
}
\]
where the horizontal maps are the duality isomorphisms, the left hand vertical map is the ordinary restriction map in cohomology and the right hand map is {\em transfer}. When $i=n$ this specializes to
\[
\xymatrix{
H^n(L,M)\ar[r]\ar[d]^{\text{Res}}&M\otimes_{\F L}\F\ar[d]^{\text{Tr}}\\
H^n(H,M)\ar[r]&M\otimes_{\F H}\F\\
}
\]
and here the transfer map is easy to describe: it is given by
$$m\otimes1\mapsto\sum_{t\in T}mt\otimes1$$ where $T$ is a left transversal to $H$ in $L$ (i.e. $L$ is the disjoint union of the left cosets $tH$, $t\in T$). We are interested in taking $M:=\F G$. Now $\F G\otimes_{\F L}\F$ is isomorphic as left $\F G$-module to the submodule of $\mathcal P$ comprising subsets which are finite unions of right cosets of $L$. Similarly $H$ is isomorphic to the left module of finite unions of cosets of $H$. From this viewpoint, the transfer map is simply induced by inclusion of sets. On passing to the colimit over all members of $\mathcal S$ we obtain Claim 2. 

We know that $H^{n+1}(G,\F G)$ is non-zero and we shall take advantage of this to construct a graph with more than one end on which $G$ acts in a useful way. We have
$$0\ne H^{n+1}(G,\F G)\iso H^1(G/\mathcal S,H^n(\mathcal S,\F G)).$$
Substituting into the five term exact sequence we obtain an exact sequence
$$0\to H^{n+1}(G,\F G)\to H^1(G,H^n(\mathcal S,\F G))\to
H^1(\mathcal S,H^n(\mathcal S,\F G))^G.$$
Thus we have an exact sequence
$$0\to H^{n+1}(G,\F G)\to H^1(G,\mathcal B)\to
H^1(\mathcal S,\mathcal B).$$
This identifies the dualizing module of $G$ over $\F$ with the kernel of a restriction map in first cohomology. To compute the first cohomology of $G$ observe that the power set $\mathcal P$ of $G$ is a coinduced $\F G$-module on which cohomology vanishes and we can view $\mathcal B$ as a submodule. The short exact sequence $\mathcal B\to\mathcal P\to \mathcal P/\mathcal B$ yields the exact sequence
$$0\to\mathcal B^G\to\mathcal P^G\to(\mathcal P/\mathcal B)^G\to H^1(G,\mathcal B)\to0$$ in cohomology. Note that $\mathcal P^G=\{\emptyset,G\}=\F$ and the assumption that all members of $\mathcal S$ have infinite index in $G$ implies that $\mathcal B^G=0$. So the exact sequence simplifies to
$$0\to\F\to(\mathcal P/\mathcal B)^G\to H^1(G,\mathcal B)\to0$$
The first cohomology group is most easily viewed as the quotient derivations modulo inner derivations. Writing $\mathcal P_{\mathcal S}$ for the preimage of $(\mathcal P/\mathcal B)^G$ under the natural map 
$\mathcal P\to\mathcal P/\mathcal B$ we have the commutative diagram
\[
\xymatrix{
&&0\ar[d]&0\ar[d]\\
&&\mathcal B\ar@{=}[r]\ar[d]&\Ider(G,\mathcal B)\ar[d]\\
0\ar[r]&\F\ar[r]\ar@{=}[d]&\mathcal P_{\mathcal S}\ar[r]\ar[d]&\Der(G,\mathcal B)\ar[r]\ar[d]&0\\
0\ar[r]&\F\ar[r]&(\mathcal P/\mathcal B)^G\ar[r]\ar[d]&H^1(G,\mathcal B)\ar[r]\ar[d]&0\\
&&0&0\\
}
\]
with exact rows and columns. The set $\mathcal P_{\mathcal S}$ consists of those subsets $B$ of $G$ such that for all $g\in G$, the symmetric difference $B+Bg$ belongs to $\mathcal B$. Such a set gives rise to a derivation defined by
$$g\mapsto B+Bg.$$ We have identified $H^{n+1}(G,\F G)$ with the kernel of the restriction map $H^1(G,\mathcal B)\to H^1(\mathcal S,\mathcal B)$. Let $\xi$ be a non-zero element of this kernel. Then $\xi$ is represented by a derivation $\delta:G\to\mathcal B$ and this restricts to an inner derivation on some subgroup $H$ in $\mathcal S$. Our derivation arises from a choice of $B\in\mathcal P$:
$$\delta g=B+Bg$$
for all $g$. The restriction condition says that there is a set $A\in\mathcal B$ such that 
$$B+Bh=A+Ah$$ for all $h\in H$. We can choose a subgroup $L$ contained in $H$ which is also a member of $\mathcal S$ such that 
$$A=LAL.$$
On restriction to $L$ we find that
$$B+B\ell=\emptyset$$ for all $\ell\in L$. This says that $B=BL$. It follows that the number of ends of the pair $G,L$ is at least $2$:
$$e(G,L)\ge2.$$
In the terminology of Dunwoody and Swenson, $L$ has codimension one in $G$ and a splitting of $G$ can be found using methods closely related to theirs, \cite{ds}.
Here we shall give an argument based on the earlier result \cite{dunwoody} of Dunwoody.

Let $X$ be a finite set of generators for $G$.
We now construct a graph $\Gamma$ with a left action of $G$. The vertex set $V$ of $\Gamma$ is the set of cosets $gL$ of $L$. 
The edge set $E$ of $\Gamma$ is defined to be a subset of $V\times V$:
$$E=\{(gL,gxL):\ g\in G, x\in X\}.$$
A typical edge $(gL,gxL)$ has initial vertex $gL$ and terminal vertex $gxL$.
Clearly $\Gamma$ admits a left action of $G$. Also $\Gamma$ is connected because $X$ generates $G$.
For each vertex $gL$ in $\Gamma$ either $gL\subseteq B$ or $gL\subseteq B^*$. 
\begin{itemize}
\item[{\bf Claim 3.}]
There are only finitely many edges $e$ of $\Gamma$ having one vertex in $B$ and one vertex in $B^*$
\end{itemize} 
To establish the claim, consider an edge $(gL,gxL)$ with $gL\subset B$ and $gxL\subset B^*$. Then $g\in B\setminus Bx^{-1}$. Similarly, if $gL\subset B^*$ and $gxL\subset B$ then $g\in Bx^{-1}\setminus B$. Thus, if the edge $(gL,gxL)$ has exactly one of its vertices in $B$ then this reasoning shows that
$$g\in \bigcup_{x\in X}B+Bx^{-1},$$ and also
$$gx\in \bigcup_{x\in X}B+Bx.$$
Set $$Y:=\left(\bigcup_{x\in X}(B+Bx^{-1})\cup (B+Bx)\right)L.$$ Then $Y$ is a union of finitely many left cosets of $L$ and every edge with exactly one vertex in $B$ has both its vertices in $Y$. This proves the Claim 3. It now follows from Dunwoody's result \cite{dunwoody} that $G$ splits over a subgroup commensurable with $L$.
\end{proof}

For Theorem C shall need to appeal to a result about group actions on trees. The following is the main theorem of \cite{bf}.

\begin{theorem}[Bestvina and Feighn (1991)]\label{bf}
Let $G$ be a group of type $\fp_2$ over $\F$. Then there exists an integer $\gamma(G)$ such that the following holds.

If $T$ is a reduced $G$-tree with small edge stabilizers, then the number of vertices in $T/G$ is bounded by $\gamma(G)$.
\end{theorem}

The precise statement in \cite{bf} assumes that $G$ is finitely presented. However, in the subsequent remark (8), the authors state that the result holds for {\em almost finitely presented groups}, i.e. groups of type $\fp_2$ over $\F$, and that their proof requires absolutely no change. The reason for this is that finite presentation is used only to manufacture a connected $2$-dimensional CW-complex $X$ on which $G$ acts freely and cocompactly in which every {\em track} separates: this last condition is guaranteed when $H^1(X,\F)=0$ and the construction of $X$ is therefore possible for any almost finitely presented $G$.

The following definitions are supplied in \cite{bf} and we restate them so that the reader can see exactly how Theorem \ref{bf} applies to our situation.

\begin{definition}\label{def:min red hyp}
Let $G$ be a group and let $T$ be a $G$-tree.
\begin{enumerate}
\item
The action of $G$ on $T$ is said to be {\em minimal} if and only if there are no proper invariant subtrees.
\item The $G$-tree $T$ is called {\em reduced} if and only if it is minimal and in addition every vertex of valency $2$ properly contains the stabilizers of the two incident edges.
\item When $T$ is minimal, it is called {\em hyperbolic} if and only if there exist two hyperbolic elements in $G$ (these being elements which have no fixed points but do have an invariant line, called the axis) whose axes intersect in a compact set. In this case, $G$ contains a free group on $2$ generators: in fact, sufficiently high powers of the two group elements freely generate a free group.
\end{enumerate}
\end{definition}

\begin{definition}\label{def:small}
A group $G$ is called {\em small} if and only if it does not admit a hyperbolic action on any minimal $G$-tree. In particular, if $G$ has no non-cyclic free subgroups then $G$ is small. Polycyclic-by-finite groups are small.
\end{definition}

The following result is a considerable generalization of the main theorem of \cite{phk-commentari}. Note that the hypothesis {\em all subgroups commensurable with $K$ are small} is clearly satisfied if $K$ is polycyclic-by-finite. The main theorem of \cite{phk-commentari} deals with the very special case when $K$ is infinite cyclic and $G$ has cohomological dimension $\le2$.

\begin{theoremC}\label{C}
Let $G$ be a group satisfying the conditions $(\alpha),(\beta),(\gamma)$ of Theorem B. Suppose that all subgroups commensurable with $K$ are small. Then there is a $G$-tree $T$ such that all vertices and edges have stabilizers commensurable with $K$.
\end{theoremC}
\begin{proof} We give the general argument below, but first, for motivation, we consider a special case which indicates how the general argument must proceed.
By Theorem B, $G$ splits over a subgroup commensurable with $K$. 
Suppose for the sake of argument that $G=U*_HV$ is a free product with amalgamation where $H$ is commensurable with $K$.
Since $G$ and $H$ are both finitely generated it is necessarily the case that $U$ and $V$ are also finitely generated. Both $U$ and $V$ satisfy all the hypotheses of Theorem $B$ and we deduce that either $|U:H|$ is finite or $U$ splits over a subgroup $L$ commensurable with $H$. If the latter holds, let $T$ be the corresponding $U$-tree. The subgroup $H$ acts on $T$ and has finite orbits on the edges of $T$. Therefore there must be at least one vertex fixed by $H$. This means that the splitting of $U$ is compatible with the original splitting of $G$ and we can find a $G$-tree combining both splittings of $G$ in which there are two orbits of edges corresponding to the subgroups $H$ and $L$. The process can be continued but we can appeal to Theorem \ref{bf} to be sure that it breaks off. When it breaks off we reach a situation where all the vertex groups as well as all the edge groups are commensurable with $K$.

More precisely, choose a reduced $G$-tree $T$ with finitely many orbits of edges, stabilizers all commensurable with $K$, and subject to these conditions, with the maximum possible number of orbits of vertices. The existence of such is guaranteed by Theorem \ref{bf}. 
This $G$-tree has finitely many orbits of edges and vertices and the fact that $G$ itself and all the edge stabilizers are finitely generated forces the vertex stabilizers to be finitely generated. 
We claim that every vertex stabilizer here is also commensurable with $K$. Suppose not. Then there is a vertex $v$ whose stabilizer $G_v$ is not commensurable with $K$. Let $e$ be an edge incident with $v$. Since $G_e$ and $K$ {\em are} commensurable we have that $G_v$ and $G_e$ are not commensurable. However, $G_e\subseteq G_v$ and therefore $G_e$ has infinite index in $G_v$. We may now apply Theorem B to split $G_v$ over a subgroup $H$ commensurable with $G_e$.
Let $T'$ be the corresponding $G_v$-tree and let $e_*$ be an edge of $T'$ with stabilizer $H$. Let $E_0$ be the set of edges which are incident with the vertex $v$ in the original tree $T$. For each $e_0\in E_0$, let $G_{e_0}$ denote the stabilizer of $e_0$ for the action of $G$ on $T$. Since $G_{e_0}\subseteq G_v$ we have an action of $G_{e_0}$ on $T'$. Moreover, $|G_{e_0}:G_{e_0}\cap H|$ is finite and so the $G_{e_0}$-orbit of $e_*$ in $T'$ is finite. It follows that $G_{e_0}$ fixes a vertex $v(e_0)$ of $T'$. We can now build a new $G$-tree 
by blowing up each of the vertices in the $G$-orbit of $v$ using the tree $T'$. 
We remove the vertices $v\cdot G$ of the orbit of $v$ and replace them with the forest $T'\times_{G_v}G$. The loose edge $e_0$ can be joined to the vertex $v(e_0)$ in the primary copy $T'\times 1$ of $T'$ and we can repeat this for the other edges incident with $v$. The process can be carried equivariantly over the orbits of loose edges. In this way we obtain a $G$-tree with a greater number of orbits of vertices. This contradicts the assumption and Theorem C follows.
\end{proof}

Other approaches to constructing the simplicial actions on trees for  Theorem B and to deducing Theorem C from Theorem B can be found in the work of Dunwoody and Swenson \cite{ds} and of Mosher, Sageev and Whyte \cite{msw1,msw2}. However the construction of almost invariant sets through the use of our new spectral sequence is novel and essential for our arguments: we do not know of any alternative to this line of reasoning beyond the two dimensional case considered in \cite{phk-commentari}.

\section{An application to Poincar\'e duality groups}

We illustrate the potential of Theorems B and C by using them to establish the following result which tidies up and extends some of the considerations in \cite{kr}. 

\begin{theorem}\label{pd}
Let $G$ be a $\PD^{n+1}$-group and let $H$ be a $\PD^n$-subgroup. Then either $|\Comm_G(H):H|$ is finite or $|G:\Comm_G(H)|$ is finite. In the latter case there is a subgroup $K$ commensurable with $H$, normal in $\Comm_G(H)$ such that $\Comm_G(H)/H$ is either infinite cyclic or infinite dihedral.
\end{theorem}

This result should not be regarded as an advance in itself. Indeed it is clear that this and further results can and have been proved by other methods in the work of Scott and Swarup \cite{ss}. However, it gives an indication of how our results may prove helpful in the study of Poincar\'e duality groups.

The following will be needed in our proof of \ref{pd}.

\begin{lemma}\label{asc}
Let $G$ be a group which is the union of a strictly ascending sequence 
$$B_0<B_1<B_2<\dots$$ of $\PD^n$-groups $B_i$. Then $G$ has cohomological dimension $n+1$.
\end{lemma}
\begin{proof}
$H^{n+1}(G,M)$ is isomorphic to $\invlim^1 H^n(B_i,M)$ for any $G$-module $M$. Taking $M=\F G$ one can calculate that this $\invlim^1 $ does not vanish.
$H^n(B_i,\F G)$ is isomorphic to $\F G/B_i$ and the connecting maps in the limit system
$$\dots\to \F G/B_i\to \dots \to \F G/B_1\to \F G/B_0$$
all strict inclusions. (See the argument for Claim 2 in the proof of Theorem B.)

We can view the inverse limit system as sitting inside the constant system
$(\F G/B_0)$ in which the connecting maps are the identity maps. We thus have a short exact sequence of limit systems
$$(\F G/B_i)\mono(\F G/B_0)\epi (U_i)$$
where the connecting maps in the quotient system $(U_i)$ are surjections with non-zero kernels. Applying the $\invlim$-$\invlim^1 $ exact sequence we obtain the exact sequence 
$$\F G/B_0\to\invlim U_i\to\invlim^1  \F G/B_i\to0.$$
The left hand map here cannot be surjective because $\F G/B_0$ is countable whereas $\invlim U_i$ has cardinality $2^{\aleph_0}$. Therefore the right hand group is non trivial, completing the proof.
\end{proof}

We shall also need Strebel's fundamental dimension theorem \cite{strebdim} which says that all subgroups of infinite index in a $\PD^{n+1}$-group have cohomological dimension $\le n$.

\begin{proof}[Proof of \ref{pd}]
As a first step we prove the
\begin{itemize}
\item[{\bf Claim.}]
Every finitely generated $S\le\Comm_G(H)$ such that $|H:H\cap S|<\infty$ is either commensurable with $H$ or of finite index in $G$.
\end{itemize}
Let $S$ be such a subgroup. We can apply Theorem B to the group $S$ with subgroup $H\cap S$. This shows that either $|S:H\cap S|$ is finite or $S$ splits over a subgroup commensurable with $H$. In the latter case, the proof of Theorem B shows that $S$ is a duality group of dimension $n+1$ over any field and so by Strebel's theorem it must have finite index in $G$.

We consider two cases which together cover all eventualities.
\begin{itemize}
\item[Case 1.]
There is a finitely generated subgroup $S$ of $\Comm_G(H)$ with 
$$|S:H|=\infty.$$
We show that in this case, $\Comm_G(H)$ has finite index in $G$ and the existence of $K$ can be established.
\end{itemize}
Applying Theorem B, we know that $S$ splits over a subgroup commensurable with $H$. 
Since the edge group in the splitting of $S$ is commensurable with $H$, it is finitely generated. Since $S$ is also finitely generated it follows that the vertex groups in the splitting of $S$ are finitely generated. The vertex groups have infinite index in $S$ and the above shows that they are commensurable with $H$. Therefore either $S=J*_KL$, a free product with amalgamation in which $J\ne K\ne L$ are all commensurable with $H$, or $S=B*_K,t$ is an HNN-extension in which the base and associated subgroups are commensurable with $H$. 

If $S$ is a non-ascending HNN-extension then one can use the Kurosh subgroup theorem to exhibit finitely generated subgroups of $S$ which have infinite index and which contain an infinite index subgroup commensurable with $H$. This contradicts the claim.

Similarly, if $G$ is an amalgamation in which one of the indices 
$|J:K|$, $|L:K|$ is $\ge3$ then we can again find intermediate finitely generated subgroups which contradict the claim.

If $S=B*_Bt$ is a strictly ascending HNN-extension then the chain of subgroups
$$B\subset B^t\subset B^{t^2}\subset\dots$$ has union of cohomological dimension $n+1$ by Lemma \ref{asc}. But this contradicts Strebel's theorem and so cannot happen. 

Therefore either $S$ is an amalgamation in which $|J:K|=|L:K|=2$ in which case $K\normal S$ and $S/K\iso D_\infty$, or $G$ is a stationary HNN-extension meaning $B=K\normal S$ and $S/K\iso C_\infty$.

Now $S$ has finite index in $G$ and normalizes $K$. Therefore $K$ has only finitely many distinct conjugates. Thus the subgroup $K_0$
defined by 
$$K_0:=\bigcap_{g\in\Comm_G(H)}K^g$$ is a finite intersection of subgroups commensurable with $H$ and is therefore itself commensurable with $H$. Also $K_0$ is normal in $\Comm_G(H)$ and the corresponding quotient is virtually cyclic. If we set $K_1/K_0$ equal to the largest finite normal subgroup of $\Comm_G(H)$ then the quotient $\Comm_G(H)/K_1$ is either infinite cyclic or infinite dihedral.

\begin{itemize}
\item[Case 2.] 
For every finite subset $F$ of $\Comm_G(H)$, 
$$|\langle H\cup F\rangle:H|<\infty.$$ 
We show that in this case, $H$ has finite index in $\Comm_G(H)$.
\end{itemize}
In this case, if $|\Comm_G(H):H|$ is infinite then we can choose a strictly ascending chain of finite extensions of $H$ by successively increasing the size of finite set $F$. The union of such a chain has cohomological dimension $n+1$ by Lemma \ref{asc} and therefore finite index in $G$. This is a contradiction because $G$ is finitely generated while the union of the chain is not. Thus
$|\Comm_G(H):H|<\infty$ as claimed.
\end{proof}

\section{An easy application in which the spectral sequence collapses}

We conclude by mentioning a very simple application of the theory motivated by the notion of complete cohomology as described in \cite{bensoncarlson,goichot,mislin}. This again is a significant generalization of one of the results in \cite{phk-commentari}.

\begin{definition}\label{def:already complete}
We shall say that a group $G$ {\em already has complete cohomology} if and only if the cohomology functors $H^n(G,{\phantom M})$ vanish on projective modules for all $n$. This is equivalent to asserting that the natural map $H^n(G,{\phantom M})\to\widehat H^n(G,{\phantom M})$ is always an isomorphism. For example, free abelian groups of infinite rank already have complete cohomology whereas finite groups never enjoy the property.
\end{definition}

\begin{theorem}\label{star}
Let $G$ be a group and let $\mathcal S$ be an admissible family of subgroups which already have complete cohomology. Then $G$ already has complete cohomology.
\end{theorem}
\begin{proof}
The new spectral sequence collapses because the hypotheses along with Lemma \ref{injectives}(iii) ensure that
$$H^*(\mathcal S,P)=0$$ for all projective modules $P$.
\end{proof}

This provides a transparent argument for proving one of the results of \cite{browngeoghegan} that
$$H^*(F,\Z F)=0$$ where $F$ denotes Thompson's group given by the presentation
$$F=\langle x_0,x_1,x_2,\dots:\ x_i^{-1}x^{{\phantom 1}}_jx^{{\phantom 1}}_i=x^{{\phantom 1}}_{j+1}\ (i<j)\rangle.$$
Let $A$ be the subgroup of $F$ generated by $\{x^{{\phantom 1}}_{2n+1}x_{2n}^{-1}:\ n\ge0\}$ and let $F_0$ be the subgroup of index $2$ in $F$ comprising elements which can be expressed as even weight words in the $x_i$. Thus
$$A\subset F_0\subset F.$$
Clearly it suffices to prove that
$$H^*(F_0,\Z F)=0,$$
and for this we only need to observe that the admissible family $\mathcal S$ of subgroups generated by $A$ consists entirely of free abelian groups of infinite rank. We shall write $A_m$ for the subgroup of $A$ generated by 
$\{x^{{\phantom 1}}_{2n+1}x_{2n}^{-1}:\ n\ge m\}$. The admissible family $\mathcal S$ is the set of subgroups which are finite intersections of conjugates $A^g$ with $g\in F_0$. Here are the precise details of the argument. 

\begin{lemma}\label{lem:Thompson}
\begin{enumerate}
\item $x^{{\phantom 1}}_1x_0^{-1}$, $x^{{\phantom 1}}_3x_2^{-1}$, $x^{{\phantom 1}}_5x_4^{-1}$, $x^{{\phantom 1}}_7x_6^{-1}$, $\dots$ is a sequence of distinct elements of $F$ which freely generate the abelian group $A$.
\item For any finite subset $X$ of $F_0$, there is an $m\ge0$ such that 
$$A_m\subseteq\bigcap_{g\in X}A^g.$$
\item Every member of $\mathcal S$ is free abelian of infinite rank.
\item $H^*(F,\Z F)=0$.
\end{enumerate}
\end{lemma}
\begin{proof}
\begin{enumerate}
\item Notice that if $0\le m<n$ then the relations in our given presentation of $F$ immediately yield
$$x_{2m}^{-1}\left(x^{{\phantom 1}}_{2n+1}x_{2n}^{-1}\right)x_{2m}
=x^{{\phantom 1}}_{2n+2}x_{2n+1}^{-1}$$
and
$$x_{2m+1}^{-1}\left(x^{{\phantom 1}}_{2n+1}x_{2n}^{-1}\right)x_{2m+1}=
x^{{\phantom 1}}_{2n+2}x_{2n+1}^{-1}.$$
Thus $x^{{\phantom 1}}_{2m+1}x_{2m}^{-1}$ commutes with $x^{{\phantom 1}}_{2n+1}x_{2n}^{-1}$. This shows that $A$ is abelian. One can see quite easily that $A$ is free abelian on the stated generators. One way is to use the known representation of $F$ as a group of piecewise linear maps of the unit interval $[0,1]$. Alternatively, let $F_m$ be the subgroup of $F$ generated by the $x_i$ with $i\ge m$. Then each $F_m$ is isomorphic to $F=F_0$ and $F_m$ is an ascending HNN-extension over $F_{m+1}$. One can now check inductively that for each $m$, the subgroup
$$x^{{\phantom 1}}_1x_0^{-1},\dots,x^{{\phantom 1}}_{2m+1}x_{2m}^{-1}$$ is free abelian of rank $m+1$ and lies in the centralizer of $F_{2m+2}$.
\item Let $g$ be a word in the alphabet $x_i^{\pm1}$, $i\ge0$ with exponent sum $j$. The relations defining $F$ show that for all sufficiently large $n$,
$$g^{-1}x_ng=x_{n+j}.$$ If $g$ is an element of $F_0$ then it can be expressed as a word with even exponent sum and therefore
$$A_m\subset A\cap A^g$$ for sufficiently large $m$. The result for an intersection of finitely many conjugates now follows from this.
\item This follows at once.
\item Theorem \ref{star} shows that $H^*(F_0,\Z F)=0$ and then we can apply the ordinary LHS spectral sequence to the group extension $F_0\mono F\epi\Z/2\Z$.
\end{enumerate}
\end{proof}

\section{Connection with Galois Cohomology}

As a concluding remark we mention that although the derived functors 
$$H^n(G/\mathcal S,{\phantom M})$$ 
are not at first sight easily related to any familiar functors, there are nevertheless close connections with Galois cohomology of profinite groups.

As an example we shall see that if $G$ is a residually finite group and $\mathcal S$ is the family of all subgroups of finite index then $\Mod\Z G/\mathcal S$ is the category of discrete modules for the profinite completion $\widehat G$ of $G$ and $H^n(G/\mathcal S,{\phantom M})$ is isomorphic to the continuous (Galois) cohomology functor $H^n_{\Gal}(\widehat G,{\phantom M})$.
We refer the reader to \cite{serrecohom} for an introduction to the theory.

In general, given a group $G$ and admissible family $\mathcal S$, 
it is not immediately clear that one can form a completion analogous to the profinite completion because it may happen that $\mathcal S$ contains few, or possibly no, normal subgroups. For example if $G$ is the Baumslag--Solitar group with presentation
$$\langle x,y:\ y^{-1}x^2y=x^3\rangle$$ and $\mathcal S$ is the set of subgroups commensurable with the infinite cyclic subgroup $\langle x\rangle$ then no member of $\mathcal S$ is normal. 

Nevertheless, 
one can always form a completion $\widehat G_{\mathcal S}$ for this example and any other. One can endow the completion with a product which makes it into a monoid. We do not know whether this monoid is necessarily a group in all cases, but we can show that it is a group in all the applications and examples considered in this paper. Define $\widehat G_{\mathcal S}$ to be the set of all functions $f:\mathcal S\to \mathcal P(G)$ which satisfy the conditions
\begin{itemize}
\item $f(H)\in H\backslash G$ for all $H\in\mathcal S$ and
\item $f(K)\subseteq f(H)$ whenever $K\subseteq H$ are members of $\mathcal S$.
\end{itemize}
In effect, we associate the coset space $H\backslash G$ to $H$ and observe that an inclusion $K\subset H$ induces a natural surjection $K\backslash G\epi H\backslash G$. Then $\widehat G_{\mathcal S}$ is the inverse limit
$\invlim\ H\backslash G$. For any $f\in\widehat G_{\mathcal S}$ and any $H\in\mathcal S$ there exists $x\in G$ such that $f(H)=Hx$. It is useful to introduce the notation $H^f$ for the conjugate $H^x$ because it depends only on $f$ and $H$; not on the particular coset representative $x$.
To make $\widehat G_{\mathcal S}$ into a monoid we define the product by setting:
$$f\cdot f'(H)=f(H)f'(H^f)$$
for each $H\in \mathcal S$. On the right, the product is carried through using the standard convention $AB=\{ab:a\in A,b\in B\}$ for sets $A,B\subseteq G$.
It is straightforward to check that $H^{f\cdot f'}=(H^f)^{f'}$.
To check the associative law:
\begin{eqnarray*}
(f\cdot f')\cdot f''(H)&=&f\cdot f'(H)f''\big(H^{f\cdot f'}\big)\\
&=&f(H)f'(H^f)f''\big(H^{f\cdot f'}\big),\\
\end{eqnarray*}
and
\begin{eqnarray*}
f\cdot(f'\cdot f'')(H)&=&f(H)(f'\cdot f'')(H^f)\\
&=&f(H)f'(H^f)f''\big((H^f)^{f'}\big).\\
\end{eqnarray*}
The function $e$ defined by $e(H)=H$ for all $H$ is the identity element.
There is a natural map $$\widehat{\phantom g}:G\to \widehat G_{\mathcal S}$$ defined by sending each $g\in G$ to the function $\widehat g$ defined by $\widehat g(H)=Hg$. This is a homomorphism which carries $1\in G$ to 
$e\in\widehat G_{\mathcal S}$ because
$$\widehat g\cdot\widehat{g'}(H)=%
\widehat g(H)\widehat{g'}(H^g)=Hgg'.$$

The question of whether or not $\widehat G_\mathcal S$ is a group appears to be subtle.
 If $f\in\widehat G_{\mathcal S}$ has an inverse $f^{-1}$ then expanding the definition of $f\cdot f^{-1}(H)$ shows that
$$f^{-1}(H^f)=f(H)^{-1}.$$
This shows that $f^{-1}$ is uniquely determined on the subset $\{H^f:H\in\mathcal S\}$ of $\mathcal S$. In general, for fixed $f$, the map
$H\mapsto H^f$ is injective and preserves the poset structure of $\mathcal S$. To see this observe that for any finite subset $\mathcal F$ of $\mathcal S$ there is a group element $x$ such that $f(H)=Hx$ and $H^f=H^x$ for all $H\in\mathcal F$. Thus the map $H\mapsto H^f$ is given locally as conjugation by a single group element.

\begin{definition}
We shall say that the admissible family $\mathcal S$ is {\em stable} if and only if for all $K\le H$ in $\mathcal S$ there exists $L\in\mathcal S$ such that $L\le K$ and $L$ is normal in $H$.
\end{definition}

\begin{lemma}
If $\mathcal S$ is a stable admissible family of subgroups of $G$ then $\widehat G_{\mathcal S}$ is a group.
\end{lemma}
\begin{proof}
Fix $f\in\widehat G_{\mathcal S}$. Fix $H\in\mathcal S$. Choose $x\in f(H)$. Choose $K$ so that $\mathcal S\ni K\le H\cap H^f$ and $K\normal H^f$. 
Choose $t\in f(K^{x^{-1}})$. Define 
$$f^{-1}(H)=Ht^{-1}.$$
\begin{itemize}
\item[{\bf Claim 1.}] $f^{-1}$ is well defined.
\end{itemize}
 To define $f^{-1}(H)$ we have made three choices, namely $x,K,t$. Note that $f(H)=Hx$, $H^f=H^x$ and $K^{x^{-1}}\le (H^f)^{x^{-1}}=H$ so necessarily $t\in f(H)$, $f(H)=Ht$ and $xt^{-1}\in H$. Thus $t^{-1}x\in H^x=H^f$ normalizes $K$ and $K^{x^{-1}}=K^{t^{-1}}$ so that $f(K^{t^{-1}})=tK$. 

Consider now a different sequence $y,L,u$ of the three choices made in the same way as $x,K,t$. Then $f(H)=Hy=Hu$ and $f(L^{u^{-1}})=uL$. Downward directedness guarantees that $tK\cap uL\ne\emptyset$. Choose $v\in tK\cap uL$. Then $tK=vK$ and $uL=vL$. Since $K$ and $L$ are contained in $H$ it follows that $Ht^{-1}=Hv^{-1}=Hu^{-1}$ as required.
\begin{itemize}
\item[{\bf Claim 2.}] $f^{-1}$ belongs to $\widehat G_{\mathcal S}$.
\end{itemize}
Given subgroups $H'\le H$ in $\mathcal S$ we can make the three choices $x,K,t$ so that $x\in f(H')$, $K\subseteq H'$ and $K\normal H$. The choices then simultaneously supply the definitions of $f^{-1}(H')$ and $f^{-1}(H)$ so that
$f^{-1}(H')=H't^{-1}$ and $f^{-1}(H)=Ht^{-1}$ have the required compatibility for the inverse system.
\begin{itemize}
\item[{\bf Claim 3.}] $f^{-1}$ is inverse to $f$. 
\end{itemize}
This follows: for any $H$ we can choose $t\in G$ such that $f(H)=Ht$, 
$f(H^{t^{-1}})=tH$, $f^{-1}(H)=Ht^{-1}$ and so
$$f^{-1}\cdot f(H)=f^{-1}(H)f(H^{f^{-1}})=Ht^{-1}tH=H=e(H).$$
Thus $f^{-1}\cdot f=e$ and since $e$ is a two-sided identity element it follows from elementary group theory that $f^{-1}$ is a two-sided inverse as required.
\end{proof}

Every object of $\Mod\Z G/\mathcal S$ can be endowed with a natural action of $\widehat G_{\mathcal S}$ as follows. 
Let $M$ be a $\Z G/\mathcal S$-module, let $f\in \widehat G_{\mathcal S}$ and let $m\in M$. Then we set
$$m\cdot f:=m.f(H)$$ where $H$ is any choice of member of $\mathcal S$ which fixes $m$. Although the subgroups of $\mathcal S$ may not be normal, it is still true that the underlying set of $\widehat G_{\mathcal S}$ can be viewed as an inverse limit:
$$\widehat G_{\mathcal S}=\limf H\backslash G$$
taken over the discrete coset spaces 
$$H\backslash G=\{Hg:\ g\in G\}.$$
Then we may endow $\widehat G_{\mathcal S}$ with the inverse limit topology and it becomes a topological group. One could then go on to consider the category of discrete $\widehat G_{\mathcal S}$-modules: this is equivalent to the category of $\Mod G/\mathcal S$-modules. The cohomology functors $H^n(G/\mathcal S,{\phantom M})$ can be identified with the continuous cohomology functors of the topological group defined for example by using continuous cocycles in the standard bar resolution construction.

Of course, if $\mathcal S$ consists of normal subgroups, or more generally if
$$\bigcap_{g\in G}H^g$$ belongs to $\mathcal S$ for all $H\in\mathcal S$, then $\widehat G_{\mathcal S}$ is easier to define: one can think straightforwardly in terms of the inverse limit of quotient groups.

 \end{document}